\documentstyle[12pt,leqno]{article}
\newtheorem{theorem}{\quad\large Theorem}
\newtheorem{proposition}[theorem]{\large Proposition}

\newtheorem{corollary}[theorem]{\large Corollary}
\newtheorem{lemma}[theorem]{\quad \large Lemma}

\input{amssym.def}

\begin{document}
\date{}

\title{Derangement Characters of the Finite General Linear Group}

\author{Alexander Gnedin\\ {\it University of Utrecht}
\and Sergey Kerov \\ {\it POMI, St. Petersburg}}
\maketitle

{\bf Abstract.} We focus on derangement characters of $GL(n,q)$ which depend solely on the dimension of the space of fixed 
vectors. This family includes  Thoma characters which become asymptotically irreducible as $n\to\infty$.
We find explicit decomposition of Thoma characters into irreducibles, construct further derangement characters
and seek for extremes in the family of derangement characters.

\vskip1cm

\par {\bf 1 Introduction.}
Let $G_n=GL(n,q)$ be the group of invertible
matrices of degree $n$ with entries from the finite field ${\Bbb F}_q$.
 Denote $1_n$ the unit matrix and for $g\in G_n$ let 
$r(g)={\rm dim\,\, ker}\, (g-1_n)$ be the number of Jordan blocks to the eigenvalue 1. 
We will consider {\it Thoma characters}
\begin{equation}\label{sigma}
\sigma_k^{(n)}(g)=q^{k\,r(g)}\qquad k=0,1,\ldots
\end{equation}
and more general {\it derangement characters}
which share with (\ref{sigma}) the
property that they depend on $g$ only through $r(g)$.

\par Our interest in the derangement characters is motivated by
the Thoma--Skudlarek  classification of  characters of the infinite
group
$GL(\infty,q)=\cup G_n$.
As was conjectured by Thoma \cite{Thoma} and then proved by Skudlarek
\cite{Skud}, any positive
definite class function $f:GL(\infty, q)\to {\Bbb C}$
satisfying $f(1_{\infty})=1$ can be uniquely
represented as a convex combination of
the functions $\sigma_k\cdot \chi$ where
$\sigma_k$ is the normalized version of (\ref{sigma}) and
$\chi$ is  a linear character.
For the 
special linear group
$SL(\infty,q)$ the one-dimensional factor is trivial and all positive definite class functions 
are of derangement type.
From a somewhat different perspective, the Thoma-Skudlarek result identifies all possible pointwise limits
of (rather complicated) characters of 
finite groups $G_n$ and suggests
that derangement characters are  the objects of
their own right.

\par 
By a {\it character} we shall mean a positive definite class function on a group.
All characters form a convex cone, which in case of finite group has traces of complex 
irreducible representations as extreme elements.
Thoma characters are  extreme as functions on the infinite linear group, but 
for fixed $n$ they can be decomposed already within the family of derangement characters.
In Section 3 we decompose Thoma characters in {\it $\psi$-characters} and then 
in Section 5  give explicit decomposition
of the $\psi$-characters into irreducibles.
The $\psi$-characters can be further decomposed within the family of derangement characters,
thus it is natural to ask about  
extreme elements of the cone of derangement characters.
One curious observation we make here is that  the cone of derangement characters is simplicial for some $n$ 
and not simplicial
for other (the first nonsimplicial case appears for $n=7$).
Although the problem of describing all extreme derangement characters for all $n$ remains open, 
we developed an algorithm
to  compute them for $n\leq 12$ (with the exception of $n=10$) and  derived a formula for  stable 
{\it $\tau-$characters}
which are extreme and span a $\lfloor n/2\rfloor$-dimensional face of the $(n+1)$-dimensional cone of 
derangement characters.

\par A straightforward generalization of Thoma characters are exponential functions $z^{n-r(g)}$.
In section 6 we give decomposition of these functions and give a new proof to  Skudlarek's result
that only values $z=q^{-k}$ yield positive definite functions for all $n$.
In Section 8 we describe a unique derangement character with trivial unipotent part.
All results of this paper  hold for arbitrary finite field.

\vspace{0.5cm}

\par {\bf 2 Characters $\sigma_k^{(n)}$.}
Let $X_k^{(n)}$ be the set of
$n$ by $k$ matrices over ${\Bbb F}_q$.
The group $G_n$ acts on $X_k^{(n)}$
by left multiplication.
For $k=1,2,\ldots$ let
$\sigma_k^{(n)}$ be the character of the
permutation representation
in ${\Bbb C}\left[X_k^{(n)}\right]$.
We also define $\sigma_0^{(n)}$ to
be the unit character of $G_n$.

\begin{proposition} The characters
$\sigma_k^{(n)}$ coincide with the Thoma characters defined by {\rm (\ref{sigma})}.
\end{proposition}
{\it Proof.} The value of $\sigma_k^{(n)}(g)$ is equal to
the number of matrices
$x\in X_k^{(n)}$ fixed by $g$.
Clearly, $x$ is fixed if and only if
each column of $x$
belongs to the space
$V_g\subset {\Bbb F}_q^{n}$ of $g$-invariant vectors.
Since the dimension of $V_g$ is $r(g)$,
there are $q^{r(g)}$ possible choices for each of $k$ columns of $x$ .$\Box$

\vspace{0.50cm}

\par Define a function of $g\in G_n$ to be {\it derangement}
if it depends on $g$ only via $r(g)$.
Obviously, any derangement function is a class 
function.
We denote $D_n$ the space of derangement functions; this
is a complex vector space of dimension $n+1$.
\par The object of our interest is  the cone $D_n^+$  of positive definite derangement
functions $f\neq 0$.
Each $f\in D_n^+$  is called a character and
${\rm dim\,} f=f(1_n)$ is called
the dimension of $f$ (even when $f$ is not a
trace of a matrix representation of $G_n$).
Because matrices $g$ and $g^{-1}$ have the same
eigenvectors we have $r(g)=r(g^{-1})$, therefore any derangement character
accepts
only real values
which are always in the range $\vert f(g)\vert \leq {\rm dim}\, f$.
\vspace{0.5cm}

\begin{proposition}
The characters $\{\sigma^{(n)}_k: 0\leq k\leq n\}$ form a linear basis
of $D_n$.
\end{proposition}
{\it Proof.} Choose arbitrary $g_j\in G_n$ with $r(g_j)=j$, $0\leq j\leq n.$
The matrix of degree $n+1$ with entries
$$(\sigma_k^{(n)})(g_j)=q^{kj}$$
is nondegenerate, because it is the Vandermonde matrix in variables
$1,q,\ldots,q^n$. It follows that
the $n+1$ characters are linearly independent
and form a basis.
$\Box$

\vspace{0.5cm}

\par We shall view $G_n$ as the generic term  of the increasing
series of groups $G_1\subset G_2\subset\ldots$ with the natural
embedding which sends $g\in G_n$ to $g\oplus 1\in G_{n+1}$.
Since the embedding adds  fixed vectors we have
\begin{equation}\label{branch-sigma}
\sigma_k^{(n)}\vert_{G_{n-1}}=q^k\,\sigma_k^{(n)}.
\end{equation}
Taken together with Proposition 2 the restriction rule implies that
the class of derangement
functions is closed under
restriction to smaller groups.

\vspace{0.5cm}

\par {\bf 3 Characters $\psi_k^{(n)}$.}
Thoma characters are not irreducible. The first obvious step to decompose them
is to split $X_k^{(n)}$ into orbits.
Let $Y_k^{(n)}$ be the space of
$n$ by $k$ matrices over ${\Bbb F}_q$ with rank $k$.
For $k=1,\ldots,n$
define $\psi_k^{(n)}$ to be the character of the permutation
representation
in ${\Bbb C\,}[Y_k^{(n)}]$
and for $k=0$ let $\psi_0^{(n)}\equiv 1$.
\vspace{0.5cm}

\begin{proposition}
The characters $\psi_k^{(n)}$ are derangement. They are given by the formula
\begin{equation}\label{psi}
\psi_k^{(n)}(g)=
(q^{r}-1)(q^{r}-q)\ldots (q^{r}-q^{k-1}), \qquad r=r(g)
\end{equation}
(which is $0$ for
$r(g)<k$).
\end{proposition}

{\it Proof.} For $V_g$ being as in Proposition 1,
a matrix $y\in Y_k^{(n)}$ is fixed by $g$ if its columns
are in $V_g$ and linearly independent.
Counting the choices compatible with the independence condition, we
have $q^r-1$ possible choices for the first column,
then $q^r-q$ choices for
the second, etc. $\Box$

\vspace{0.5cm}

\par The following branching rule is  analogous to the restriction
formula
(\ref{branch-sigma}).

\begin{proposition} For $1\leq k\leq n$ we have
\begin{equation}\label{branch-psi}
\psi_k^{(n)}\vert_{G_{n-1}}= q^k \psi_k^{(n-1)}
+(q^{2k-1}-q^{k-1})\psi_{k-1}^{n-1}
\end{equation}
(the first term is void for $k=n$).
\end{proposition}
{\it Proof.} Fix $g\in G_{n-1}$ with $r(g)=g$. The embedding
$G_{n-1}\subset G_n$ sends
$g$ to $g\oplus 1$ with $r(g\oplus 1)= r+1$.
We shall count the number of matrices $x\in Y_k^{(n)}$ fixed
by $g\oplus 1$. Let $x$ be such a matrix, and ${\hat x}$ be this
matrix with the last row deleted, then of course $g{\hat x}={\hat x}$
and the rank of ${\hat x}$ must be either $k$ or $k-1$.
\par In the first case there are $\psi_k^{(n-1)}(g)$ choices for
${\hat x}$ which could be arbitrarily combined with any of $q^k$
choices for the last row of $x$. This yields the first term in
(\ref{branch-psi}).

\par In the second case ${\hat x}$ can be
seen as a $k$-tuple of column vectors
which span a $(k-1)$-dimensional subspace
of $V_g$ (${\rm dim}\, V_g=r$). Making further distinction
between the cases
when the first column of ${\hat x}$ belongs to the space spanned by the rest
$k-1$ columns or not we compute the number of choices for ${\hat x}$ as
$$(q^r-1)\ldots (q^r-q^{k-2})\frac{q^k-1}{q-1},$$
which is a multiple of $\psi_{k-1}^{(n-1)}(g)$ by Proposition 3.
Since the rows of ${\hat x}$ also span a $(k-1)$-dimensional space,
there are always $q^k-q^{k-1}$ ways to extend ${\hat x}$ to a matrix of rank
$k$ by appending the last row. This results in the second term
in (\ref{branch-psi}). $\Box$

\vspace{0.5cm}

\par Recall that the $q$-binomial coefficient is defined as
$${k\choose j}_q=\frac{(q^k-1)(q^{k-1}-1)\ldots (q^{k-j+1}-1)}
{(q^{j}-1)(q^{j-1}-1)\ldots (q-1)}$$
and is equal to the number of $j$-dimensional subspaces in ${\Bbb F}_q^k$
(which is 0 for $k<j$).
\vspace{0.5cm}

\begin{proposition} The characters $\sigma_k^{(n)}$ and
$\psi_k^{(n)}$  are related by the formulas
\begin{equation}\label{sigmapsi}
\sigma_k^{(n)}=\sum_{j=0}^k {k\choose j}_q \psi_j^{(n)}
\end{equation}
\begin{equation}\label{psisigma}
\psi_k^{(n)}=\sum_{j=0}^k (-1)^{k-j}q^{{k-j\choose 2}}
{k\choose j}_q
\sigma_j^{(n)}
\end{equation}
\end{proposition}
{\it Proof.}  We  decompose $X_k^{(n)}$ into $G_n$-orbits.
Let $x\in X_k^{(n)}$ be a matrix of rank $j$ with columns
seen as elements of ${\Bbb F}_q^n$.
Select $j$ linearly independent columns,
label them  $v_1,\ldots,v_j$ and label $v_{j+1},\ldots,v_k$ the rest
columns.
Consider the space of linear relations $\alpha$ in $k$ indeterminates
over ${\Bbb F}_q$ such that $\alpha (v_1,\ldots,v_k)=0$
and let $\alpha_1,\ldots,\alpha_{n-j}$ be a basis of this space.
Since the natural action of $G_n$ in ${\Bbb F}_q^n$ is $j$-transitive
(for $j\leq n$) the orbit $\{gx: g\in G_n\}$ coincides with
the set of matrices
whose columns satisfy the $\alpha_i$'s and the first $j$ columns are
independent.
It is easily seen that the orbit of $x$ is isomorhic to
$Y_k^{(n)}$.
Using the isomorphism
$X_k^{(n)}\approx {\Bbb F}_q^n\otimes {\Bbb F}_q^k$,
the linear span of the orbit becomes
${\Bbb F}_q^n\otimes A$
where $A\subset {\Bbb F}_q^{k}$ is the null space for
relations $\alpha_i$.

\par Any
choice of independent columns results in the same
$A$, therefore
the correspondence between the orbits and subspaces of ${\Bbb F}_q^k$
is bijective. It follows that
$X_k^{(n)}$ splits into
${k \choose j}_q$ orbits isomorphic to $Y_j^{(n)}$, $j=0,\ldots,n$,
which implies
(\ref{sigmapsi})
(the term $j=0$ corresponds to the singleton
orbit of the zero matrix in $X_k^{(n)}$).
The second formula (\ref{psisigma}) follows from the first by virtue of
an inversion formula for the $q$-binomial coefficients
(analogous to a better
known inversion formula for the binomial coefficients).  $\Box$

\vspace{0.5cm}
\par It follows that the characters $\{\psi_k^{(n)}: 0\leq k\leq n\}$
also form a linear basis
of $D_n$ and that
both sets of characters generate the same integral lattice.

\par The equivalence of branching
formulas (\ref{branch-sigma}) and (\ref{branch-psi}) can be also
derived from the expansions (\ref{sigmapsi}) and (\ref{psisigma}).
In fact, the relation between
the two sets of characters
is the specialization of
the $q$-binomial formula
$$(1+\xi)(1+\xi q)\ldots (1+\xi^{k-1}q)=\sum_{j=0}^k {k \choose j}_q
q^{{j\choose 2}} \xi^j$$
for $\xi=-q^{r-k+1}$.

\vspace{0.5cm}

\par {\bf 4 On irreducible characters of $G_n$.}
We will need  some well-known facts about the
irreducible characters of $G_n$,
referring the reader to \cite{Ze}, \cite{MD} for  a fuller account.

\par As usual, we identify
Young diagram
$\lambda=(\lambda_1,\ldots ,\lambda_m)$
with its geometric image
$\{(i,j):1\leq i\leq m, 1\leq j \leq \lambda_i\}$,
and write
$\vert\lambda\vert=\lambda_1+\ldots +\lambda_m$
for the number of boxes.
We denote ${\Bbb Y}$ the collection
of all Young diagrams,  including  the empty diagram with
$\vert\emptyset\vert= 0$.

\par Given integers $1\leq k<n$  and two characters
$f_1$ and $f_2$ of the groups $G_k$ and $G_{n-k}$, respectively,
their parabolic product $f_1\circ f_2$ is
the character of $G_n$
induced from the parabolic subgroup
\begin{equation}
\label{parab}
P=\left\{\left({g_{1}\atop 0}{* \atop
g_{2}}\right):  g_{1}\in G_k,\, g_{2}\in G_{n-k}\right\}
\end{equation}
by the function $f_1(g_{1})\cdot f_2(g_{2})$
(in this context the groups are embedded in
$G_n$ as suggested by the definition of $P$).
The dimension of the parabolic product is
\begin{equation}\label{dimprod} {\rm dim}\, f_1\circ f_2={n\choose k}_q {\rm
dim}\, f_1 \cdot {\rm dim}\, f_2 \end{equation}
as it follows from the
Frobenius formula for induced characters and the observation that the left
coset classes for $P$ can be labeled by $k$-dimensional subspaces $V\subset
{\Bbb F}_q^n$.

\par A character is called cuspidal if it is not a part of any
parabolic product. We denote ${\cal C}_d$ the finite set of cuspidal
characters of $G_d$ and ${\cal C}=\cup_{d\geq 1} {\cal C}_d$.
The unit character
of $G_1={\Bbb F}_q^{\times}$
plays a distinguished role and will be denoted $e$;  it
is one of $q-1$
elements of ${\cal C}_1$.
Given a family $\varphi:{\cal C}\to {\Bbb Y}$  with
finitely many nonvoid diagrams $\phi(c)\neq\emptyset$, its degree is defined as
$$\| \varphi \|=\sum_{d\geq1 }\sum_{c\in {\cal C}_d} d\cdot\vert \varphi
(c)\vert.$$ A fundamental fact says that the irreducible characters of $G_n$
are in one-to-one correspondence with the families of Young diagrams of degree
$n$.  We denote $(\varphi)$ the character corresponding to such a family
$\varphi$.

\par The character of $G_n$
corresponding to the family with a single nonvoid diagram
$\phi(e)=\lambda$ is called unipotent and will be denoted
$(\lambda)_e$.
For unipotent characters the dimension
can be computed by the
$q$-hook formula:

\begin{equation}\label{qhook}
\dim\, (\lambda)_e=\frac{(q-1)\ldots
(q^n-1)}{\prod_{b\in\lambda}(q^{h(b)}-1)} \,\, q^{n(\lambda)} \,\,,
\qquad\qquad
n=\vert\lambda\vert \end{equation}
where $h(b)$ denotes
the hook length
$\lambda_i+\lambda_j'-i-j+1$ at box $b=(i,j)$ and $n(\lambda)=\sum
(i-1)\lambda_i$.
The involved function of $q$ is also called
the Kostka-Foulkes polynomial $K_{\lambda, (1^n)}(q)$.

\vspace{0.5cm}
\par {\large Remark.} To avoid confusion
between the diagram and its transpose
keep in mind that we adopt the
parametrization which relates the one-row diagram to the
unit character $(n)_e$,
while the one-column diagram corresponds
to Steinberg character $(1^n)_e$
of dimension $K_{(1^n),(1^n)}(q)=q^{{n\choose 2}}$.
In many sources the convention is reverse
(e.g. \cite{MD}).

\vspace{0.5cm}

\par {\bf 5 Decomposition of $\psi_k^{(n)}$ into blocks of irreducibles.}
Let ${\rm reg}_k$ be the character of the regular representation
of $G_k$. The next proposition says that, in a sense,
characters $\psi_k^{(n)}$ interpolate between
the unit and the regular character of $G_n$.
\begin{proposition}
We have $\psi_0^{(n)}=(n)_e\,,
 \,\,\psi_n^{(n)}={\rm reg}_n\,,$ and for
$1\leq k\leq n$
\begin{eqnarray}
\psi_k^{(n)}=(n-k)_e\circ{\rm reg}_k\,.
\end{eqnarray}
\end{proposition}
{\it Proof.} Since $G_n$ acts transitively on $Y_k^{(n)}$ character
$\psi_k^{(n)}$ is induced
by the unit character from the
subgroup $B$ of block matrices
$\left({g_1\atop 0}{0\atop 1_k} \right), g_1\in G_{n-k},$
which fix $\left({0_{(n-k)\times k}\atop 1_k}\right)\in Y_k^{(n)}$.
The parabolic product structure is recognized
when we view the induction as the two-step procedure:
at first inducing from $B$ to the parabolic group
(\ref{parab}) - which yields the  character
$(n-k)_e (g_1) \cdot {\rm reg}_k(g_2)$ of $P$ -
and then using this character to further
induce from $P$ to $G_n.$ $\Box$

\vspace{0.5cm}

\par Introduce the character
$\rho_j=\sum_{{\| \varphi \|=j\atop \varphi(e)=\emptyset}}
 \dim (\varphi) \cdot (\phi)$
of $G_j$ and
for $\lambda\in {\Bbb Y},\vert\lambda\vert\leq n$ set
\begin{equation}\label{lambdablock}
[\lambda]_n =(\lambda)_e\circ \rho_{n-\vert\lambda\vert}
\end{equation}
Decomposition of $[\lambda]_n$ into irreducibles follows directly
from the definition: this involves all family of Young diagrams with
$\phi(e)=\lambda$.

\par Characters $[\lambda]_n$
are disjoint in the sense that none of the irreducible characters
of $G_n$ enters decompositions of two different $[\lambda]_n$'s.
We shall see that $[\lambda]_n$'s are convenient `blocks'
for  representing  derangement characters,
they are themselve irreducible if $\lambda$ has $n$ boxes
(in which case $[\lambda]_n=(\lambda)_e$).

\par Given $\mu\in {\Bbb Y}$ and integer $m$ let ${\rm H}_m^+(\mu)$
be the set of Young diagrams which can be derived  from $\mu$
by appending
a horizontal strip.
That is to say, $\lambda\in {\rm H}_m^+(\mu)$
if the skew diagram $\lambda\setminus \mu$ has $m$ boxes with
no two in the same column.  Reciprocally, let
${\rm H}_m^-(\lambda)$ be the set
of diagrams which can be derived  from $\lambda$ by
deleting a horizontal strip.
Clearly,
\begin{equation}\label{HH} \lambda\in
{\rm H}_m^+(\mu)\quad \Longleftrightarrow\quad \mu\in {\rm H}_m^-(\lambda).
\end{equation}
Note that the largest horizontal strip has $\lambda_1$ boxes,
thus ${\rm H}_m^-(\lambda)$ is empty if the first row of $\lambda$
is shorter than $m$.
\par Next result gives explicit decomposition into blocks.

\begin{theorem} Each character $\psi_k^{(n)}, 0\leq k\leq n,$
is an integral linear combination of the characters
$[\lambda]_n$:
$$\psi_k^{(n)}=\sum_{\vert\lambda\vert\leq n}
c_k^{(n)}[\lambda]_n.
$$
The multiplicity is zero if  $\vert\lambda\vert>n$ or
$\lambda_1<n-k$,
otherwise
\begin{equation}\label{c}
c_k^{(n)}(\lambda)={k\choose n-\vert\lambda\vert}_q\,\,
\sum_{\mu\in {\rm H}_{n-k}^- (\lambda)}
\dim\, (\mu)_e \end{equation}
\end{theorem}
{\it Proof.} The regular character ${\rm reg}_k$
splits into irreducibles as
\begin{eqnarray*}\label{reg}
{\rm reg}_k=
\sum_{\| \varphi \|=k}
\dim (\varphi)\cdot (\varphi)=\\
\sum_{\vert\mu\vert\leq k}
\sum_{\varphi (e)=\emptyset\atop \| \varphi \|= k-\vert\mu\vert}
\dim\, ((\mu)_e\circ
(\varphi))\cdot (\mu)_e\circ (\varphi)=\\
\sum_{\vert\mu\vert \leq k}
{k\choose \vert\mu\vert}_q
\dim\,(\mu)_e\cdot
(\mu)_e\circ\rho_{k-\vert\mu\vert}
\end{eqnarray*}
as it follows from (\ref{dimprod}) and (\ref{lambdablock}).

\par  For $\nu\in {\Bbb Y}$ and integer $m$ the parabolic
product decomposes as
\begin{equation}\label{prod}
(m)_e\circ (\nu)_e=\sum_{\lambda\in {\rm H}^+_m(\nu)} (\lambda)_e.
\end{equation}
This follows by virtue of Pieri's rule
which is the same for 
unipotent characters of $G_n$
 as for 
characters of the symmetric group.
From this, Proposition 6 and (\ref{lambdablock})
we obtain

\begin{eqnarray*}
\psi_k^{(n)}={\rm reg}_k\circ (n-k)_e=\\
\sum_{\vert\mu\vert\leq k} {k\choose\vert \mu\vert}_q
\dim\, (\mu)_e
\sum_{\lambda\in {\rm H}_{n-k}^+(\mu)}
\rho_{k-\vert\mu\vert}\circ (\lambda)_e=\\
\sum_{\vert\mu\vert\leq k} {k\choose\vert \mu\vert}_q
\dim (\mu)_e
\sum_{\lambda\in {\rm H}_{n-k}^+(\mu)} [\lambda]_n.
\end{eqnarray*}

The coefficient at $[\lambda]_n$ is calculated by swapping the sums,
applying (\ref{HH}) and observing that
$n-\vert\lambda\vert=k-\vert\mu\vert$ implies
$${k\choose \vert\mu\vert}_q={k\choose n-\vert\lambda\vert}_q.$$
$\Box$

\vspace{0.5cm}

\par  In one most important case
the multiplicity formula (\ref{c}) simplifies. Given
$\lambda$ let $\nu=\lambda\setminus (\lambda_1)$ be the diagram
derived from $\lambda$ by
deleting the first row. Note that deleting the maximum horizontal
strip of $\lambda$ also yields $\nu$.
Suppose $i=\lambda_1-(n-k)\geq 0$, then removing a
horizontal strip with $n-k$ boxes from $\lambda$
is equivalent to appending
a horizontal strip with $i$ boxes to $\nu$.
On the other hand, appending $i$ boxes to $\nu$ results in some diagram
$\mu\in {\rm H}_{n-k}^-(\lambda)$ provided that $\mu_1\leq \lambda_1$.
It follows that (\ref{c}) is equivalent to
$$c_k^{(n)}(\lambda) ={n\choose n-\lambda}_q
\sum_{\mu\in {\rm H}_i^+(\nu)\atop \mu_1\leq\lambda_1} \dim\,(\mu)_e. $$

\begin{corollary} If $\lambda_2\leq n-k\leq \lambda_1$
then
$$c_k^{(n)}(\lambda)={k\choose j}_q
{k-j\choose i}_q\dim \nu$$
where $\nu=\lambda\setminus (\lambda_1)$, $i=\lambda_1-(n-k)$  and
$j=n-\vert\lambda\vert$.
If $n\geq 2k$ then for any $\lambda$ with at most $n$ boxes
\begin{equation}\label{csimple}
\psi_k^{(n)}=\sum_{{i+j\leq k\atop i,j\geq
0}} {k \choose j}_q {k-j\choose i}_q\,\,
\sum_{\vert\nu\vert=k-i-j} \dim \nu
\cdot [n-k+i,\nu].
\end{equation}
\end{corollary}
{\it Proof.} Suppose $\lambda_2\leq n-k\leq \lambda_1$ then
$\lambda_2+i=\lambda_2+\lambda_1-(n-k)\leq\lambda_1$ and therefore
$\mu_1\leq \lambda_1$ for any $\mu\in {\rm H}_i^+(\nu)$.
By Pieri's rule and (\ref{dimprod}) we get
$$\sum_{\mu\in {\rm H}_i^+(\lambda)}\dim\, (\mu)_e=\dim\, (\nu)_e\circ (i)_e=
{k-j \choose i}_q\dim\, (\nu)_e$$
because $\vert\nu\vert+i=\vert\lambda\vert -\lambda_1 +i=(n-j)-(i+n-k)+i=k-j$.
\par If $n\geq 2k$ the inequality $\lambda_1\geq n-k$
implies $\lambda_2\leq n-\lambda_1\leq n-(n-k)=k\leq n-k$.
Now (\ref{csimple}) follows because
any diagram entering the decomposition of $\psi_k^{(n)}$ is of the form
$\lambda=(n-k+i,\nu)$ with $\vert\nu\vert\leq k-i$.$\Box$
\vspace{0.5cm}

\par {\large Remark.} The coefficient (\ref{c}) is a
multiple of the skew Kostka polynomial\\
$K_{\lambda\setminus (n-k),(1^{k-j})}(q)$,
as introduced in \cite{Kir}.
Under conditions of Corollary 8 the
skew diagram $\lambda\setminus (n-k)$ splits in
two parts with no common boxes in the same row
or column; factoring of the polynomial also follows
from the interpretation as the generating function of
tableaux (see \cite{Kir}).
\vskip0.5cm
\par A positivity property of Kostka-Foulkes polynomials implies that
(\ref{c}) are polynomials with positive integral coefficients.

\begin{corollary}

\begin{itemize}
\item[{\rm (i)}] Character $[\lambda]_n$ with the first row $\lambda_1=n-k$ is present only in the decomposition of
$\psi_k^{(n)},\ldots ,\psi_n^{(n)}$.
\item[{\rm (ii)}] The empty diagram enters
only the regular character, so that
$c_k^{(n)}(\emptyset)=\delta_{kn}$ (Kronecker
delta).
\item[{\rm (iii)}] We have for the regular character
$$\psi_n^{(n)}=
c_n^{(n)}(\lambda)={n\choose\vert\lambda\vert}_q\dim\,(\lambda)_e
[\lambda]_e.$$
\item[{\rm (iv)}] For one-row diagrams we have
$c_k^{(n)}(n-k)=1$, and 
more generally for `hook diagrams'
$\lambda=(n-k+i,1^{k-i-j})$ with $i,j\geq 0;\,i+j\leq k<n$:
$$c_k^{(n)}(\lambda)={k\choose i}_q{k-i\choose j}_q q^{k-i-j\choose 2}.$$


\end{itemize} \end{corollary}
{\it Proof.} Straightforward from
(\ref{c}).  For (iv) apply (\ref{csimple}).$\Box$
\vspace{0.5cm}

\par {\large Example.} We tabulate coefficients of the  decomposition
into $[\lambda]_n$'s for $n=4$.
\[
\begin{array}{lccccc}
      & \psi^{(4)}_0& \psi^{(4)}_1& \psi^{(4)}_2& \psi^{(4)}_3& \psi^{(4)}_4\\

(4)   &   1         &      1      &       1     &       1     &    1        \\

(3,1) &    0         &      1      &      1+q    &  1+q+q^2    &    q+q^2+q^3 \\

(2^2) &   0          &  0          &      1      &    q+q^2    &    q^2+q^4  \\

(2,1^2)& 0          &   0& q& q+q^2+q^3& q^3+q^4+q^5\\

(1^4) &   0&0&0&q^3&q^6\\

(3)  & 0&1&1+q&1+q+q^2&1+q+q^2+q^3\\

(2,1)&0&0&1+q&1+2q+2q^2+q^3& q+2q^2+2q^3+2q^4+q^5\\

(1^3)& 0&0&0& q+q^2+q^3&q^3+q^4+q^5+q^6\\

(2)&  0&0&1&1+q+q^2&1+q+2q^2+q^3+q^4\\

(1^2)& 0&0&0& 1+q+q^2& q+q^2+2q^3+q^4+q^5\\

(1) & 0&0&0& 1& 1+q+q^2+q^3\\

\emptyset &0&0&0&0&1
\end{array}
\]

\vspace{0.5cm}

\par Any derangement function $f$ has a unique representation as
a linear combination of the $\psi_k^{(n)}$'s and this
implies a decomposition of $f$ into $[\lambda]_n$'s.
In particular, Theorem 7 combined with Proposition 5
enables representing Thoma characters
$\sigma_k^{(n)}$ as integral linear combination of the $[\lambda]_n$'s.

\vspace{0.5cm}

\par {\bf 6 Generalized Thoma characters.} Given $z\in {\Bbb C}$
consider the derangement function $f_z(g)=z^{n-r(g)}$.
This definition is consistent for different $n$
because $r(g\oplus 1)=r(g)+1$, thus
$f_z$ is defined on the infinite group $G_{\infty}=\cup G_n$.
Obviously, for $z=0$,
$f_z$ is the normalized
regular character equal to $\delta_{1_{\infty},g}\,$, while
for $z=q^{-k}$ it is the normalized
Thoma character
$$\sigma_k:=\frac{\sigma_k^{(n)}}{\dim\,\sigma_k^{(n)}}=q^{-k(n-r(g))}.$$

\par Skudlarek proved that for $z\neq 0$ the only positive definite
functions among $f_z$ are Thoma characters (see \cite{Skud},
Behauptung 3). His proof exploited an embedding of the
additive group of infinite matrices into $G_{\infty}$.
We show next that this result follows rather easily from the
decomposition of $f_z$ into irreducible characters for $n=1,2,\ldots$.

\begin{proposition} For $z\in {\Bbb C}$
\begin{equation}\label{f_z}
f_z=\sum_{j=0}^n\left(
z^{n-j}\prod_{i=0}^{j-1}\frac{q^{-i}-z}{q^{i+1}-1}
\right) \psi_j^{(n)}.
\end{equation}
\end{proposition}
{\it Proof.} By (\ref{sigmapsi})
$$\sigma_k=q^{-kn}\sum_{j=0}^n {k\choose j}_q \psi_j^{(n)}$$
for $k=0,\ldots,n$, which transforms into
(\ref{f_z}) with $z=q^{-k}$.
But this implies that (\ref{f_z}) holds everywhere because
for each $r(g)$ both parts of the formula are polynomials in $z$.$\Box$

\begin{corollary} The function $f_z:G_{\infty}\to {\Bbb C}$
is positive definite if and only if $z=0$ or $z=q^{-k},\,\,k=0,1,\ldots$.
\end{corollary}
{\it Proof.} If $f_z$ is positive definite then for each $n$
the coefficients in the decomposition into $[\lambda]_n$'s
are nonnegative. Character $[\emptyset]_n$
enters only $\psi_n^{(n)}$ with the coefficient being a positive
multiple of $(1-z)(q^{-1}-z)\ldots (q^{-n+1}-z),$ which in turn
is positive
for all $n$ provided that
either $z$ is from the conjectured list or $z<0$.
Hence we only need to exclude negative values.
\par Steinberg character $[1^n]_n$ enters only $\psi_{n-1}^{(n)}$ and
$\psi_n^{(n)}$.
By (\ref{f_z}) and Corollary 9 (iv)
the coefficient at $[1^n]_n$ is
$$
q^{{n\choose 2}}\prod_{i=0}^{n-1}\frac{q^{-i}-z}{q^{i+1}-1}+
q^{{n-1 \choose 2}}z \prod_{i=0}^{n-2}\frac{q^{-i}-z}{q^{i+1}-1}\,.$$
For $z<0$ positivity of the coefficient amounts to the inequality
$$z\geq \frac{-1}{q^n-q^{n-1}-1}$$
which has the right-hand side vanishing as $n\to\infty$ ($q>1$); hence
$f_z$ cannot be positive definite on all $G_n$'s if $z<0$.$\Box$

\vspace{0.5cm}

\par {\bf 7 The cone $D_n^+$.} All characters of $G_n$ form a cone
whose extreme rays correspond to irreducible characters.
This cone is simplicial, so that
any character has a unique representation as a positive
linear combination of the irreducibles.
Since derangement characters are always reducible
(besides $(n)_e$)
it is natural to ask which of them are `the least reducible'.

\par A derangement character $f\in D_n^+$ is said to be {\it extreme}
if $f',f-f'\in D_n^+$ implies $f=pf'$ for some $p>0$;
in this case $\{pf:p>0\}$ is called the extreme ray.
(Sometimes we make no distinction between a character
and its positive multiples.)
Extreme characters are of primary interest because each
$f\in D_n^+$ is a positive linear combination
of the extremes.

\par Let $L_n$ be the space of
class functions representable as
linear combinations of the $[\lambda]_n$'s and $L_n^+$ be the cone
of characters in $L_n$. 

\par Although $[\lambda]_n$'s are not derangement (see next proposition), they offer a useful coordinatization of
the space of derangement functions.
Indeed,
because $L_n$ is generated by disjoint
irreducible characters each $f\in L_n$ has a unique representation
as
$$f=\sum_{\vert\lambda\vert\leq n}f\langle\lambda\rangle [\lambda]_n$$
and $f\in L_n^+$ if and only if the coefficients $f\langle\lambda\rangle$
are nonnegative.
We will call the set of diagrams
${\rm supp\,}f:=\{\lambda: f\langle\lambda\rangle
>0\}$ {\it the support} of
$f\in L_n^+$.

\par Because $D_n\subset L_n$  we have
$D_n^+=D_n\cap L_n^+$ and we can use this fact
to distinguish the characters from other derangement functions.
The cone $D_n^+$ is polyhedral and has a compact base (a polytope).
Since $D_n^+$
contains linearly
independent characters $\psi_k^{(n)},0\leq k\leq n,$
the number of extreme rays must be at
least $n+1$, and if it is exactly $n+1$
the cone is simplicial and each character has a unique
representation as a positive linear combination of extremes.

\par By Corollary 9 (i) we see that 
$${\rm supp\,}\psi_k^{(n)}=\{\lambda:\lambda_1\geq n-k\},$$
and this implies an important observation that 
{\it the supports are strictly increasing with $k$}.
An immediate consequence is
\vskip0.5cm

\begin{proposition}
None of the characters $[\lambda]_n$ is derangement,
besides the unit character $[n]_n$.
\end{proposition}
{\it Proof.} Given $f\in D_n$, let $n-j$ be the length of the shortest
first row of all $\lambda$'s entering $f$ with some nonzero coefficient.
Since the supports are increasing, $j$ is the maximum
index of the nonzero $a_k$'s entering $f=\sum a_k\psi_k^{(n)}$.
But then all $[\lambda]_n$ with $\lambda_1=n-j$ enter $f$ with the same
$a_j$. Now, if $f=[\lambda]_n$ then $\lambda_1$ is the
shortest first row of all diagrams.  But for $0<\lambda_1<n$
there are other diagrams with the same first row
which would enter in the derangement case.
That $[\emptyset]_n\not\in D_n$ will follow from Theorem 17.
$\Box$
\vspace{0.5cm}

\par Let $C_k$ be the cone obtained via intersecting  the space spanned
by  $\psi_0^{(n)},\ldots ,\psi_k^{(n)}$ with $D_n^+$.
The cones
$C_k$ are increasing with $k$,  and from the increasing of supports follows
that each cone $C_k$ is a $(k+1)-$dimensional
face of $D_n^+$, whence the following claim.

\begin{lemma} $f\in C_k$ is extreme in $D_n^+$ if and only if
$f$ is an extreme character in $C_k$.
\end{lemma}

\par Given a finite set $\{f_j:j\in J\}\subset L_n$ and $i\in J$ we say that
$\lambda$ is an {\it eigendiagram} of $f_i$ if
$\lambda\in {\rm supp\,}f_i$
but $\lambda \not\in \cup_{j\in J\setminus \{i\}}\, {\rm supp\,}f_j$.

\begin{lemma}
$C_k$ is simplicial if and only if there is
a list of characters $f_0,\ldots,f_k\in C_k$
such that  each $f_i$ has an eigendiagram (in which case this is
the complete list of extremes).

\end{lemma}
{\it Proof.} Suppose each $f_i$ has an eigendiagram,
then the linear mapping which assigns
$f_0,\ldots,f_k$
to the basis vectors of ${\Bbb R}^{k+1}$
is an order isomorhism sending $C_k$ to the positive orthant.
It follows that there are no other extreme characters.
\par Conversely, suppose $C_k$ is simplicial and let $f_0,\ldots,f_k$
be a complete list of extreme elements such
that $f_0$ has no eigendiagram.
In this case ${\rm supp\,}f_0\subset \cup_{i\neq 0}\,{\rm supp\,}f_j$
thus selecting $a_i> 0$ sufficiently large we obtain some
$f=\sum_{j\neq 0}a_j f_j-f_0\in C_k$.
We get then  for
$f_0+f\in C_k$ one decomposition $\sum_{i\neq 0} a_if_i$ without
$f_0$, while
decomposing $f$ we get another decomposition
$f_0+f=f_0+\sum b_if_i$ (with $b_i\geq 0$) which does involve $f_0$.
Since for simplicial cone the decomposition into extremes
must be unique we have a contradiction.
$\Box$

\vspace{0.5cm}

\par There is a simple method to  verify if a character
is extreme.

\begin{lemma}
Let $f_0\ldots, f_n$ be a basis of $D_n$ with $f_0\in D_n^+$.
The character $f_0$ is extreme if and only if the column vectors
$$\{f_j\langle\lambda\rangle :\lambda\not\in {\rm supp\,}f_0\}\quad
j=1,\ldots,n$$
are linearly independent. Hence each extremal character must have at least
$n$ zero coefficients.
\end{lemma}
{\it Proof.} If linear independence does not hold there is a
linear combination $f=\sum_{i=1}^n a_jf_j$ such that
$f\langle\lambda\rangle=0$ for all $\lambda\in {\rm supp\,}f_0$.
Selecting $\epsilon>0$ we obtain a character $f=f_0-\epsilon f\in D_n^+$
with ${\rm supp\,}f'\in {\rm supp\,}f_0$.
For $p$ sufficiently large $pf_0-f'\in D_n^+$.
Thus $f_0$ is not extreme since the characters $f_0$ and $f'$ are not
colinear.

\par If the linear independence does hold there is no noncolinear
character $f$ with ${\rm supp\,}f \subset {\rm supp\,}f_0$, because
there is no linear combination as above.$\Box$
\vspace{0.5cm}

\par These considerations motivate introducing yet another basis
$\tau_0^{(n)},\ldots,\tau_n^{(n)}\in D_n$
which we define recursively, as the
output of the following elimination algorithm.

\vspace{0.5cm}
\par {\large The elimination algorithm.} Set
$\tau_0^{(n)}:=\psi_0^{(n)}$. At each stage $k=1,\ldots ,n$
we have characters $\tau_j^{(n)},j<k,$ at hand and determine
sequentially characters $\tau_{k0},\ldots,\tau_{kk}$ by setting at first
$\tau_{k0}:=\psi_k^{(n)}$ and for $j=1,\ldots ,k-1$
$$\tau_{kj}:=\tau_{k,j-1}-a_{k,j-1}\tau_{j-1}^{(n)}$$
where the coefficient $a_{k,j-1}$ takes the maximum possible value compatible
with the condition that the difference be a character (i.e. in $L_n^+$).
Finally,  define $\tau_k^{(n)}:=\tau_{kk}.$
\vspace{0.5cm}

\par {\large Remark.} Explicitly,
the coefficients
are
$$a_{i,j}=\min_{\lambda}
\tau_{k,j}\langle\lambda\rangle /\tau_{j-1}^{(n)}\langle\lambda\rangle$$
($a/0=\infty$).
Complemented by $a_{i,j}=\delta_{ij}, i\leq j,$ they determine the
transition matrix from the basis $\{\tau_k^{(n)}\}$ to
$\{\psi_k^{(n)}\}$.
\par To apply the algorithm one needs
to determine the minima like $\min P_{\lambda}(q)$
for certain polynomials in $q$.
However, there is a computer evidence that this problem is trivial:
the polynomials involved have positive coefficients and there is
always a polynomial which has minimal coefficients
at all powers of $q$.

\vskip0.5cm

\begin{proposition}
If the cone $C_k$ is simplicial then $\tau_0^{(n)},\ldots, \tau_k^{(n)}$
is the complete list of extreme characters of $C_k$.
\end{proposition}
{\it Proof.} The statement is trivial for $k=0$.
If $C_k$ is simplicial then the same applies to $C_{k-1}$, so suppose
by induction that $\tau_0^{(n)},\ldots, \tau_{k-1}^{(n)}$
are extreme in $C_{k-1}$ (thus, by obvious extension of Lemma 13, are extreme also in $C_k$).
Note that $C_k$ is in the linear span of
$\psi_k^{(n)},\tau_0^{(n)},\ldots, \tau_{k}^{(n)}$.
Consider the step resulting in $\tau_{k1}$.
In geometric terms, the elimination means that we determine
the intersection point
of the ray connecting $\tau_0^{(n)}$ and $\psi_k^{(n)}$ with
the face $Q\subset C_k$ not containing $\tau_0^{(n)}$.
Obviously, $Q$ is a simplicial cone of lower dimension,
the characters $\tau_1^{(n)},\ldots,\tau_k^{(n)}$
are extreme in $Q$ and the linear span of
$\tau_{k1},\tau_1^{(n)},\ldots,\tau_k^{(n)}$ contains $Q$.
This is the same situation as with $C_k$ but now the dimension
is reduced, thus the induction step can be completed. $\Box$
\vspace{0.5cm}

\par Combining Lemma 14 and the last proposition, we see that $D_n^+$ is
simplicial provided each $\tau_k^{(n)}$ has an eigendiagram,
otherwise
not. Next example illustrates the approach.

\vskip0.5cm

\par {\large Example.} The decomposition of characters $\tau_k^{(4)}$ ($n=4$) is

     \[ \begin{array}{lccccc} & \tau^{(4)}_0& \tau^{(4)}_1&
      \tau^{(4)}_2& \tau^{(4)}_3& \tau^{(4)}_4\\

(4)   &   1         &      0      &       0     &       0     &    0        \\

(3,1) &    0         &      1      &      0      &   0         &   0          \\

(2^2) &   0          &  0          &      1      &    0        &   0         \\

(2,1^2)& 0          &   0& q& q        & 0          \\

(1^4) &   0&0&0&q^3&0  \\

(3)  & 0&1&0  &0      &1          \\

(2,1)&0&0&1+q&1+ q         & q+ q^2              \\

(1^3)& 0&0&0& q+q^2+q^3&q^3            \\

(2)&  0&0&1&1      &1+q+ q^2        \\

(1^2)& 0&0&0& 1+q+q^2&   q  + q^2+q^3    \\

(1) & 0&0&0& 1& 1+q+q^2    \\

\emptyset &0&0&0&0&1
\end{array}
\]

It is seen that $(4),(3,1),(2^2),(1^2),\emptyset$ are the eigendiagrams.
Therefore $D_4^+$ is simplicial, and $\tau_k^{(4)},$ $0\leq k\leq 4,$ is the
complete list of extreme derangement characters.
The linear relations between the bases are the following:
\[
\begin{array}{lrrrrr}
\psi_0^{(4)}=& \tau_0^{(4)} &&&&\\
 \psi_1^{(4)}=& \tau_0^{(4)} &+   \tau_1^{(4)}&&&\\
\psi_2^{(4)}=& \tau_0^{(4)} &+  (q+1) \tau_1^{(4)}&+  \tau_2^{(4)}&&\\
\psi_3^{(4)}=&\tau_0^{(4)} &+  (q^2+q+1) \tau_1^{(4)}
&+ (q^2+q)\tau_2^{(4)}& + \tau_3^{(4)}&\\
\psi_4^{(4)}=&\tau_0^{(4)} &+  (q^3+q^2+q+1)
\tau_1^{(4)} &+  (q^4+q^2)\tau_2^{(4)}&+ q^3\tau_3^{(4)}&+ \tau_4^{(4)}.
\end{array}
\]

\vspace{0.5cm}

\par Results of similar computations for $n\leq 22$ are as follows:

\vskip0.5cm

\begin{itemize}
{\it
\item[{\rm (i)}] For $n=1,2,3,4,5,6,8,9,11,12$ the
cone $D_n^+$ is simplicial
and $\{\tau_k^{(n)}, 0\leq k\leq n\}$ is the complete list of extreme
characters.

\item[{\rm (ii)}] For $n=7, 10,13,14,15,16,17,18,19,20,21, 22$
the cone is not simplicial.
All characters $\tau_k^{(n)}$ are extreme but this list is not complete.

\item[{\rm (iii)}] For $n=7,10,13,14,15,16,17,18,19,20$,
the number of characters $\tau_k^{(n)}$ with no
eigendiagram equals one, for $n=21$ this number is two and for $n=22$ it is three.  

}
\end{itemize}

\vskip0.5cm
\par In all above cases if
 $\tau_k^{(n)},k<n,$ has an eigendiagram it is of almost rectangular shape $((n-k)^a,b)$ where $n=a(n-k)+b$
(i.e. the shape differs from rectangular only in the last row).

\vspace{0.5cm}
\par {\large Example.} The minimum $n$ such that  $D_n^+$ is not
simplicial is $n=7$.
There are 9 extreme derangement characters:
$\tau_k^{(7)}, 0\leq k\leq 7,$ and one additional character
$$\tau_*^{(7)}:=
a_1\tau_4^{(7)} +a_2 \tau_6^{(7)} -\tau_5^{(7)},$$
where
$$a_1= \frac{(1+q)(1+q^2)}{(1+q+q^2)}\,\,\,\,\,{\rm and}\,\,\,\,
a_2=\frac{(1+q)(1+q^2)^2(1+q+q^2+q^3+q^4)}{q^2+q^4+q^5+q^6+q^7+q^8+q^{10}}.
$$
The table shows the pattern of positive coefficients at
unipotent $[\lambda]_n$ (with diagrams of full degree $n$)
and $[\emptyset]_n$

\[
\begin{array}{lccccccccccc} {\rm sign}      & \tau^{(7)}_0& \tau^{(7)}_1&
\tau^{(7)}_2& \tau^{(7)}_3& \tau^{(7)}_4   & \tau^{(7)}_5 &
\tau^{(7)}_6 &
\tau^{(7)}_7 &&
\tau^{(7)}_*                \\

(7)      &+&0&0&0&0&0&0&0  && 0              \\
(6,1)    &0&+&0&0&0&0&0&0  && 0   \\
(5,2)    &0&0&+&0&0&0&0&0  && 0       \\
(5,1^2)  &0&0&+&0&0&0&+&0  && +         \\
(4,3)    &0&0&0&+&0&0&0&0  && 0        \\
(4,2,1)  &0&0&0&+&0&+&+&0  && +        \\
(4,1^3)  &0&0&0&+&0&+&+&0  && +          \\
(3^2,1)  &0&0&0&0&+&0&0&0  && +           \\
(3,2^2)  &0&0&0&0&+&+&0&0  && 0      \\
(3,2,1^2)&0&0&0&0&+&+&+&0  && +       \\
(3,1^4)  &0&0&0&0&+&+&+&0  && +       \\
(2^3,1)  &0&0&0&0&0&+&+&0  && 0       \\
(2^2,1^3)  &0&0&0&0&0&+&+&0  && +         \\
(2,1^5)  &0&0&0&0&0&+&+&0  && +           \\
(1^7)    &0&0&0&0&0&0&+&0  && +        \\
\emptyset &0&0&0&0&0&0&0&+ && 0        \\

\end{array}
\]
The eigendiagrams are $(7),(6,1),(5,2),(4,3),(3^2,1),(1^7),\emptyset$,
but $\tau_5^{(7)}$ has no eigendiagram.
Each row in the rest of the coefficients matrix
is a positive linear
combination of the rows of the above block.
The completeness of the list
of extreme characters was shown
with the help of Lemma 15.
The base of the $D_7^+$ is a polytope which is
combinatorially equivalent to a
5-fold pyramid build upon the `square'
$\{\tau_4^{(7)},\tau_5^{(7)}, \tau_*^{(7)},\tau_6^{(7)}\}$.

\vspace{0.5cm}

\par Our computations strongly suggest the following

\vskip0.5cm

\par {\large Unipotent conjecture:} each  row  of the coefficients matrix corresponding to a diagram
with less than $n$ boxes  is a positive linear combination of
the rows corresponding to unipotent characters and to $[\emptyset]_n$
(for $n\leq 22$ it was sufficient to take unipotent characters with almost rectangular shape).
An equivalent property is that the cone dual to
$D_n^+$ is  spanned by positive combinations of the rows of the unipotent block.
Equivalently, the projection $\tau\mapsto
\sum_{\vert\lambda\vert=n} \tau\langle\lambda\rangle [\lambda]_n$
is an isomorphism of ordered spaces.

\vspace{0.5cm}

\par {\bf 8 The character with no unipotent part.} Introduce
the character
\begin{equation}\label{nonuni}
{\hat \tau}_n^{(n)}=\sum_{\vert\lambda\vert\leq
n-1}{n-1\choose \vert\lambda\vert}_q \dim\, (\lambda)_e\cdot [\lambda]_n\,.
\end{equation}
In all cases covered by the computational results of previous section,
this character coincides with $\tau_n^{(n)}$.
Although we failed to prove that the coincidence is not incidental
we will show that ${\hat \tau}_n^{(n)}$ is indeed derangement and give
it characterization.

\par Denote $n_q=1-q^n$ and $n_q!=1_q 2_q\cdots n_q\,$.

\begin{theorem} Character ${\hat \tau}_n^{(n)}$ is
extreme. It can be characterized as the unique
(up to a scalar multiple) derangement function orthogonal to all
unipotent characters $[\lambda]_n$.
In terms of the basis
derangement characters:
\begin{equation}\label{decomp}
{\hat \tau}_n^{(n)}=
\psi_n^{(n)}-\sum_{k=0}^{n-1}
\frac{(n-1)_q!}{k_q!}\,q^k \psi_k^{(n)}\,.
\end{equation}
\end{theorem}

\par Note that the sum in (\ref{decomp}) is alternating, since $k_q<0$
for $q>1$. The proof of this result is based on one nontrivial
identity.

\begin{lemma} For each diagram $\lambda$ with $n$ boxes
\begin{equation}\label{identity}
c_n^{(n)}(\lambda)-
\sum_{k=0}^{n-1}
\frac{(n-1)_q!}{k_q!}\,q^k c_k^{(n)}(\lambda ) = 0
\end{equation}
\end{lemma}
{\it Proof.} Schur functions satisfy
\begin{equation}\label{schur}
s_{\lambda}(1,\xi_1,\xi_2,\ldots)=\sum s_{\mu}(\xi_1,\xi_2,\ldots)
\end{equation}
where the summation is over the set of diagrams
${\rm H}^-(\lambda)=\{\mu:\lambda\setminus (\lambda_1)\subset
\mu\subset\lambda\}$ which can be derived from $\lambda$ by deleting a
horizontal strip (the term with $\mu=\lambda$ is also
included in the right-hand side).
Similar formula with $m$ variables amounts to the
branching rule for characters of $GL(m,{\Bbb C})$.
Specializing the
Schur function for $\xi_j=q^j$ we have
$$s_{\lambda}(1,q,q^2,\ldots)=\frac{\dim\,(\lambda)_e}{n_q!}$$
and by homogeneity $s_{\mu}(q,q^2,\ldots)=q^{\mu}s_{\mu}(1,q,q^2,\ldots)$
(see \cite{Stanley}, p. 375).

\par In view of Theorem 7 the left-hand side of (\ref{schur}) is
$c_n^{(n)}(\lambda)/n_q!$, while the right-hand side is

$$\sum_{k=0}^n \sum_{\mu\in {\rm H}_{n-k}^-(\lambda)}
c_n^{(k)}(\lambda)\frac{q^k}{k_q!}  $$
which implied readily (\ref{identity}).$\Box$
\vskip0.5cm
\par {\large Remark.}
Formula (\ref{identity}) is a hidden version
of the
Kostka-Foulkes polynomials identity found in \cite{Kir2}, with a minor
correction.  Our proof is borrowed from \cite{Kir2}
(where the formula needs correction by  taking $\lambda$ in place of transpose $\lambda'$).

\vspace{0.5cm}
\par {\large Example.} For hook diagrams $\lambda=(n-m,1^m)$
the identity amounts to
$${n-1\choose m}_q\, q^{{m+1\choose 2}}=
\sum_{k=0}^{n-1}
\frac{(n-1)_q!}{k_q!}\,q^k\, {k\choose m}_q\, q^{{m\choose 2}}$$
Simplifying this becomes
$$ \sum_{j=0}^N q^j \,\frac{N_q!}{j_q!}=1$$
which can be proved straightforwardly by induction on $N$.
\vspace{0.5cm}

{\it Proof of Theorem 17.} For $\lambda$ with less than $n$ boxes
set $j=n- \vert\lambda\vert $ and observe the recurrence

\begin{equation}\label{c-reduce}
c^{(n)}_k(\lambda)={k\choose j}_q c^{(n-j)}_{k-j}(\lambda)
\end{equation}
(which is $0$ for $k<j$). Plugging this into the left-hand side of
(\ref{identity}) and applying Lemma 18 for $n'=n-j$ along with the
$q-$binomial
identity
${n\choose j}_j={n-1\choose j-1}_q+q^j{n-1\choose j}_q$
and $c_{n-j}^{(n-j)}=\dim\, (\lambda)_e$ we obtain
$$c_n^{(n)}(\lambda)-
\sum_{k=0}^{n-1}
\frac{(n-1)_q!}{k_q!}q^k c_k^{(n)}(\lambda)=
{n-1\choose \vert\lambda\vert}_q \dim\,\lambda$$
which is equivalent to (\ref{decomp}) by (\ref{nonuni}) and Theorem 7.
\par It is obvious from the definition and (\ref{decomp}) that
${\hat \tau}_n^{(n)}\in D_n^+$.
For each $k=1,\ldots, n-1$ the  character $\psi_k^{(n)}$ includes
at least one unipotent $[\lambda]_n$ which is not in
 $\cap_{j<k}{\rm supp\,}\psi_{j}^{(n)}$.
Hence the linear
rank of the `unipotent' matrix block $U$ with
entries
$$U_{k,\lambda}=c_k^{(n)}(\lambda)\qquad \vert\lambda\vert=n,\,\,
0\leq k\leq n$$
cannot exceed $n$.
The uniqueness claim now follows because
there are no two noncolinear combinations
of the basis derangement characters with zero unipotent part.
We see that the rank of $U$ is $n$, hence in
accord with Lemma 15 ${\hat\tau}_n^{(n)}$ is extreme.
$\Box$

\vspace{0.5cm}

\par As a by-product we obtain:

\begin{corollary} For $j=0,1,\ldots,n$, the linear rank of the matrix block
$U(n-j)=\{c_k^{(n)}(\lambda): \vert\lambda \vert=n-j, 0\leq k\leq n\}$
is $n-j$.
\end{corollary}
{\it Proof.} We have proved this for the unipotent block;
and for other blocks this follows from
(\ref{c-reduce}) by induction on $n$. $\Box$

\vskip0.5cm

\par It is not at all obvious from the explicit formula
$${\hat \tau}_n^{(n)}(g)= (-1)^nq^{{n\choose 2}}n_q! \,\delta_{rn}\,  -  (n-1)_q! r_q!\,,\qquad g\in G_n,\, r=r(g)$$
that this function is positive definite.

\vspace{0.5cm}

\par {\large Remark.}
It would be interesting to learn if there is some
intrinsic relation
between ${\hat \tau}_n^{(n)}(g)$  and regular characters.
The complement in ${\rm reg}_n$
is yet another derangement character
$${\rm reg}_n-{\hat \tau}_n^{(n)}=
\sum_{\vert\lambda\vert\geq 1}  q^{n-\vert\lambda\vert}
{n-1\choose
\vert\lambda\vert-1}_q\
\dim\, (\lambda)_e \,[\lambda]_n,$$
but {\it formal} replacing $[\lambda]_n$ by $[\lambda]_{n-1}$
in (\ref{nonuni}) yields ${\rm reg}_{n-1}$.
\vspace{0.5cm}

\par {\bf 9 Stable characters.} In Section 7 we introduced characters $\tau_k^{(n)}$ implicitly, by 
a recursive procedure. For $k\leq n/2$ there is an explicit formula
\begin{equation}\label{stabletau}
\tau_k^{(n)}=\sum_{j=0}^k {k\choose j}_q\sum_{\vert\mu\vert=j} \dim\, (\mu)_e
\cdot[n-k,\mu]_n\,
\end{equation}
and the relation of these characters with $\psi_k^{(n)}$'s is the same $q-$binomial as the relation between
$\psi-$  and $\sigma -$characters (\ref{sigmapsi}):
\begin{equation}\label{psitau}
\psi_k^{(n)}=\sum_{j=0}^k {k\choose j}_q \tau_j^{(n)}
\end{equation}
as it  follows easily from the simple case (\ref{csimple})
of the formula for coefficients.

\par To prove that these characters indeed appear as the output of the algorithm, we 
can just start by {\it defining} them by one of the two formulas. Then we observe that $\tau_k^{(n)}$
has eigendiagram $(n-k,k)$ (which other characters of this set do not have), therefore 
they span $C_{\lfloor n/2\rfloor}$, they are extreme by Lemma 14 and this cone is simplicial. By Proposition 16,
the elimination algorithm gives the full list of  extremes in the simplicial case, thus these extreme characters
coincide with (\ref{stabletau}) (possibly up to a positive factor).

\par From Lemma 13 follows that
\vskip0.5cm
\begin{theorem}
The characters $\tau_k^{(n)}$ are extreme for $k\leq n/2$.
\end{theorem}
\vskip0.5cm

\par We wish to stress that (\ref{psitau}) is only valid for the indicated range, and inverting the formula
for $k\leq n/2$ would not produce positive definite functions at all.
\par Inverting (\ref{psitau}) and applying (\ref{branch-psi})
we get yet another  branching rule
\begin{equation}\label{branch-tau}
\tau_k^{(n)}\vert_{G_{n-1}}= q^k \tau_k^{(n-1)}
+2q^{k-1}(q^{k}-1)\tau_{k-1}^{n-1}
+q^{k-2}(q^{k-1}-1)(q^k-1)\tau_{k-2}^{n-1}
\end{equation}
(with obvious adjustments for extreme values of indices).

\par Characters (\ref{stabletau}) are {\it stable} in the sense that, as $n$ grows, the diagrams 
entering  the decomposition of such a character keep changing only in the number of boxes in the first row.
Asymptotic considerations, which lie outside the scope of this paper, show that the normalized characters
$\psi_k^{(n)},\tau_k^{(n)}$ approach Thoma characters as $n\to\infty$.

e-mail: gnedin@math.uu.nl

\end {document}